\documentclass{article}
\usepackage{todonotes}
\author{Luis Fernando Garc\'ia-Mora and Hugo Alberto Rinc\'on Mej\'ia.}
\title{Co-first $\rmod$.}
\date{}

\usepackage{amsmath}
\usepackage{footnote}
\usepackage{multicol}
\usepackage{booktabs}
\usepackage[flushleft]{threeparttable}
\usepackage[font=small,labelfont=bf]{caption}
\usepackage{todonotes}

\usepackage[utf8]{inputenc}
\usepackage{multicol}

\usepackage{blindtext}
\usepackage{amssymb}
\usepackage{amsthm}
\usepackage{amsmath}
\usepackage{dsfont}
\usepackage{mathrsfs}
\usepackage{tikz-cd}
\usetikzlibrary{cd}
\usepackage{comment}

\newcommand{\tfor}{\mathbb{T}}

\newcommand{\trad}{\textbf{$R$-trad}}
\newcommand{\trid}{\textbf{$R$-trid}}
\newcommand{\rpr}{\textbf{$R$-pr}}
\newcommand{\zpr}{\textbf{$\mathbb{Z}$-pr}}

\newcommand{\rmod}{R\textit{-Mod}}
\newcommand{\rsimp}{R\textit{-simp}}
\newcommand{\rid}{R\textit{-pid}}
\newcommand{\rrid}{R\textit{-rid}}
\newcommand{\rrad}{R\textit{-rad}}

\newcommand{\rlep}{R\textit{-lep}}

\newcommand{\rtors}{R\textit{-tors}}
\newcommand{\rTORS}{R\textit{-TORS}}
\newcommand{\rconat}{R\textit{-conat}}

\newcommand{\Mmid}{\hspace{.7mm}\vert\hspace{.7mm}}

\theoremstyle{plain}
\newtheorem{theorem}{Theorem}[section]
\newtheorem{corollary}[theorem]{Corollary}
\newtheorem{lemma}[theorem]{Lemma}

\newtheorem{proposition}[theorem]{Proposition}

\theoremstyle{definition}

\newtheorem{definition}[theorem]{Definition}
\newtheorem{example}[theorem]{Example}

\newtheorem{remark}[theorem]{Remark}

\begin{document}
\begin{center}
{\Large \bf  Co-first modules.} 

\vspace{6mm}

{\large \bf Luis Fernando Garc\'ia-Mora and Hugo Alberto Rinc\'on-Mej\'ia.}

\vspace{3mm}
\end{center}

\begin{multicols}{2}

\begin{flushright}\noindent\textbf{MSC2020:}\hspace{1cm}\phantom{.}\\
16S90,16S99,16K99 \end{flushright}

\end{multicols}
\begin{abstract}

This paper explores the concept of \textbf{co-first modules}, a generalization of coprime modules, through the lens of preradicals in module theory.  We introduce new submodule products and characterize coprime modules using these products.  The study extends classical definitions by defining \textbf{$\mathscr{A}$-co-first modules} and \textbf{$\mathscr{A}$-fully co-first modules}, which utilize subclasses of preradicals to broaden the scope of coprimeness.  We investigate the lattice structure of conatural classes and their closure properties, providing conditions under which co-first modules coincide with second modules.  We examine the implications of these for left MAX rings and left perfect rings, comparing coprimeness and secondness.

\end{abstract}
\section*{Introduction}
In recent years, several concepts related to primality for $R$-modules have been introduced and studied, providing valuable insight into various structural properties. Of particular interest is the \emph{dual notion of prime submodule}, such as the concept of a \emph{second module} introduced by S.~Yassemi in \cite{Yas}. A nonzero submodule $N$ of $_R M$ is called \emph{second} if for every ideal $I$ of $R$, either $IN = 0$ or $IN = N$. This notion was later extended via classes of preradicals in \cite{Sec}.\\

A pretorsion class of modules $\zeta$ is a class of modules closed under  quotients and coproducts. A pretorsion class $\zeta$ is a torsion class if it is closed under extensions. A torsion class $\zeta$ is an hereditary torsion class if it is closed under submodules. It is well know that the class of hereditary torsion classes is a a lattice which for the purposes of this paper we will denote as $\rtors$. The class of torsion classes forms a big lattice  which we will denote as $\rTORS$.\\

Another important concept is that of \emph{coprime modules}, introduced by L.~Bican, T.~Kepka, and P.~N\v{e}mec in Section 4 of \cite{Bic1}. A module $M$ is said to be \emph{coprime} if every nonzero quotient of $_R M$ generates the pretorsion class $\xi(M)$, where $\xi(M)$ denotes the (not necessarily hereditary) pretorsion class generated by $M$. Explicitly,
\[
\xi(M) = \left\{ K \;\middle\vert\; \text{there exists an epimorphism } g : M^{(X)} \to K \text{ for some set } X \right\}.
\]

In the same section of \cite{Bic1}, the authors introduce the following product of two submodules $A$ and $B$ of a given module $M$:
\[
A \Box_M B = \bigcap \left\{ f^{-1}(A) \;\middle\vert\; f \in \mathrm{Hom}_R(M,M) \text{ and } f(B) = 0 \right\}.
\]
It is shown that a module $M$ is coprime if and only if for all proper submodules $L, N$ of $_R M$, the product $L \Box_M N$ is a proper submodule of $M$ (see Proposition~4.3 in \cite{Bic1}).

In this paper, we propose a generalization of the submodule product using preradicals, investigate its properties, and show how coprime modules can be characterized through this product. Furthermore, we introduce the concept of \emph{$\mathscr{A}$-co-first modules}, which generalizes the notion of coprime modules by means of a generalization of both the submodule product proposed in \cite{Bic1} and the one introduced here.

 In \textbf{Section~1}, we present some background concepts and results on preradicals that will be useful throughout the paper. 
 
In \textbf{Section~2}, we observe that the submodule product described in Section~4 of \cite{Bic1} can be expressed in terms of certain preradicals. We also introduce a new submodule product and use it to characterize coprime modules. 

In \textbf{Section~3}, we define the notion of \emph{$\mathscr{A}$-coprime modules}, which extends the classical notion of coprimeness via classes of preradicals, and compare it with the standard definition. We also define \emph{$\mathscr{A}$-fully co-first modules}. 

In \textbf{Section~4}, we study the classes of \emph{$\sigma$-fully co-first modules}, which turn out to be conatural classes, and explore several related properties. 

Finally, in \textbf{Section~5}, we compare the extension of the notion of coprime modules developed in this article with the extension of second modules proposed in \cite{Sec}, and analyze conditions on the ring under which these two notions coincide.

\section{Preradicals in $\rmod$.}

In this section, we introduce the foundational concepts of preradicals and strongly invariant submodules, as well as some results that will prove useful throughout the remainder of the article.

We introduce the basic definitions and results of preradicals in $\rmod$. For more information on preradicals, see \cite{Bic}, \cite{Clark}, \cite{pre}, \cite{preII}, \cite{preIII}, and \cite{Es}.\begin{definition}
A preradical $\sigma$ on $\rmod$ is a functor $\sigma:\rmod\rightarrow \rmod$ such that:
\begin{enumerate}
    \item $\sigma(M)\leq M$ for each $M\in \rmod$.
    \item For each $R$-morphism $f:M\rightarrow N$,  the following diagram is commutative: \begin{center}
\begin{tikzcd}
 M \ar{r}{f}& N\\
\sigma(M)\ar[hookrightarrow]{u}\ar{r}{f_\downharpoonright^\upharpoonright} & \sigma(N).\ar[hookrightarrow]{u}
\end{tikzcd}
\end{center}
\end{enumerate}
\end{definition}
Recall that for each $\beta\in \rpr$ and each $\{M_i\}_{i\in I}$ family of $R$-modules, we have that $\beta (\bigoplus \limits _{i\in I}M_i)=\bigoplus \limits _{i\in I}\beta (M_i)$, see Proposition I.1.2 of \cite{Bic}. 
$\rpr$ denotes the collection of all preradicals on $\rmod$. By Theorem 7 of \cite{pre}, $\rpr$ is a big lattice with the smallest element $\underline{0}$ and the largest element $\underline{1}$.
    In  $\rpr$, the order is defined by \( \alpha \preceq \beta \) if \( \alpha(M) \leq \beta(M) \) for every \( M \in \rmod\). Furthermore, \( (\bigvee_{i \in I} \sigma_i)(M) = \sum_{i \in I} \sigma_i(M) \) and \( (\bigwedge_{i \in I} \sigma_i)(M) = \bigcap_{i \in I} \sigma_i(M) \).\\
    
Let $\sigma \in \rpr$. A module $M$ is a $\sigma$ pretorsion module if it satisfies $\sigma(M) = M$, while it is a $\sigma$ pretorsion free module if $\sigma(M) = 0$. The collection of all $\sigma$ pretorsion modules is denoted by $\mathbb{T}_\sigma$, and the collection of all $\sigma$ pretorsion free modules is denoted by $\mathbb{F}_\sigma$.\\

 Recall that $\sigma\in \rpr$ is idempotent if $\sigma \circ \sigma =\sigma$.  $\sigma$ is radical if $\sigma(M/\sigma(M))=0$, for each $M\in \rmod$.  $\sigma$ is a left exact preradical if it is a left exact functor. 
$\sigma$ is $t$-radical (also known as a cohereditary preradical) if $\sigma(M)=\sigma(R)M$. 
Recall that $\sigma$ is a $t$-radical if and only if $\sigma$ preserves epimorphisms, see Exercise 5 of Chapter VI \cite{Es}.     $\sigma$  is a left exact preradical if and only if, for each submodule $N$ of a module $M$, we have $\sigma(N)=\sigma(M)\cap N$; see Proposition 1.7 of Chapter VI of \cite{Es}. \\
\textbf{Subclasses of preradicals}

We will denote by $\rid$, $\rrad$, $\rlep$, $\trad$, and $\trid$, the collections of idempotent preradicals, radicals, left exact preradicals, $t$-radicals, and idempotent $t$-radicals, respectively. We denote by $\rsimp$ a set or representatives of isomorphism classes of simple $R$-modules.\\
For $\sigma \in \rpr$, we denote by $\widehat{\sigma}$ the largest idempotent preradical smaller than or equal to $\sigma$ and $\overline{\sigma}$ the least radical greater or equal than $\sigma$, following Stenström. See  \cite{Es}, Chapter VI, Proposition 1.5.\\
\textbf{Fully invariant modules.}

Recall that, a submodule $N$ of $M$ is a fully invariant submodule of $M$, if $f(N)\subseteq N$ for each $f\in End_R(M)$. The set of fully invariant submodules of a module $M$ forms a complete sublattice of $\mathscr{L}(M)$, where $\mathscr{L}(M)$ is the lattice of submodule of $M$.

For a fully invariant submodule $_RN$ of $_RM$, the preradicals $\alpha_N^M$ and $\omega_N^M$ are defined in \cite{pre}, Definition 4.  
\begin{definition}
Let $N$ be a fully invariant submodule of $M$ and $U\in \rmod$. 
\begin{enumerate}
\item $\alpha_N^M(U)=\sum\lbrace f(N) \Mmid f\in Hom_R(M,U) \rbrace$ and 
\item $\omega_N^M(U)=\bigcap\lbrace f^{-1}(N) \Mmid f\in Hom_R(U,M) \rbrace$.
\end{enumerate}  
\end{definition}

\begin{remark}
Let $M\in \rmod$ and let $N$ be a fully invariant submodule of $M$. By \cite{pre} Proposition 5,  $\alpha_N^M :\rmod\rightarrow \rmod$ is the least preradical $\rho$ such that $\rho(M)=N$, and $\omega_N^M :\rmod\rightarrow \rmod$ is the largest preradical $\rho$  such that $\rho(M)=N$. The set $\{\rho\in \rpr\mid \rho(M)=N\}$ is the interval $[\alpha_N^M, \omega_N^M]$ in $R$-pr.
\end{remark}
\begin{remark}\label{tryrej}
The preradical $\alpha_M^M(L)$ is idempotent and is the same thing as the trace of $M$ in $L$, specifically $\alpha_M^M = tr_{M}$. Furthermore, $\omega^M_0(L)$ coincides with the rejection of $M$ in $L$ which is defined as the smallest module of the set $\{K \leq L \mid L/K \text{ embeds in } M^X \text{ for some set } X\}$.
Since $N/\omega_0^M(N)$ embeds in a product of copies of $M$, the least submodule of $N/\omega_0^M (N$ with this property is $0$. Thus $(\omega_0^M:\omega_0^M)(N)/\omega_0^M(N)=\omega_0^M(N)/\omega_0^M(N)$. Hence  $(\omega_0^M:\omega_0^M)=\omega_0^M$, therefore, $\omega_0^M$ is a radical. 
\end{remark}

\section{Comultiplication, totalizers and coprime modules}

Recall that for a submodule $N$ of $M$, we have the preradical 
$\gamma_N^M\colon \rmod \to \rmod$, 
where for every $U \in \rmod$ the following hold:
\[
\gamma_N^M(U) = 
\bigcap\left\{
  f^{-1}(N) \mid 
  f \in \mathrm{Hom}_R(U,M)
\right\},
\]

With this in mind, observe that for any submodules $A$ and $B$ of $M$, 
we have that $A\Box_MB$ is the unique submodule of $M$ such that 
\[
(A \Box_M B)/B =\gamma_A^M(M/B).
\]
We now consider a new product obtained by using the radical  
$\omega_0^{M/A}$ in place of the preradical $\gamma_A^M$, 
as described in the following definition:
\begin{definition}
     Let $A,B$ be submodules of $_RM$. We define $A:B$ as s  the unique submodule of $M$ such that:\begin{center}
         $(A:B)/B:=\omega^{M/A}_0 (M/B)$.
     \end{center}
     We will call $A:B $ the comultiplication of $A$ and $B$.
\end{definition} 
 Note that $(M:N)/N=M/N$, since $(M:N)/N=\omega^{M/M}_0 (M/N)=M/N$  because $\omega^0_0$  is the  identity in $R$-mod. 

The following example shows that the comultiplication is different from the product $\_\Box_M\_$.

\begin{example}
 Let us consider the submodules $2\mathbb{Z}$ and $6\mathbb{Z}$ of $\mathbb{Z}$. Since $\mathbb{Z}_2$ is embedded in $\mathbb{Z}_6$ we have that $\omega_0^{\mathbb{Z}_6}(\mathbb{Z}_2)=0$. On the other hand, since $\mathbb{Z}$ has no torsion, there are no non-zero morphisms of  $\mathbb{Z}_2$ to $\mathbb{Z}$, thus  $\gamma_{6\mathbb{Z}}^\mathbb{Z}(\mathbb{Z}_2)=\mathbb{Z}_2$. Therefore, $6\mathbb{Z}\Box_\mathbb{Z}2\mathbb{Z}=\mathbb{Z}$ and $6\mathbb{Z}: 2\mathbb{Z}=2\mathbb{Z}$.
\end{example}
\begin{proposition}
    Let $A,B,C\in \mathscr{L}(M)$. If $B\leq C$, then $(A:B)\leq (A:C)$.
\end{proposition}

\begin{proof}
    We consider the $R$-morphism $\pi:M/B\rightarrow M/C$ given by $\pi(x+B)=x+C$. Then the following diagram is commutative:
    \begin{center}
    \begin{tikzcd}
 M/B \ar{r}{\pi}& M/C\\
(A:B)/B=\omega_0^{M/A}(M/B)\ar[hookrightarrow]{u}\ar{r}{\pi_\downharpoonright^\upharpoonright} & \omega_0^{M/A}(M/C)=(A:C)/C.\ar[hookrightarrow]{u}
\end{tikzcd}
\end{center}
    Then $(A:B)+C\leq (A:C)$, so $(A:B)\leq (A:C)$.
\end{proof}

\begin{corollary}
    Let $A,B\in \mathcal{L}(M)$ and  $\{N_i\}_{i\in I}\subseteq \mathcal{L}(M)$.  Then \begin{center}
        $(A:\bigcap\limits _{i\in I}N_i)\leq \bigcap\limits _{i\in I}(A:N_i)$.
    \end{center}
\end{corollary}

\begin{lemma}\label{tot}
       There is a smallest submodule $_RU$ of $_RM$  such that $U:N=M$.
\end{lemma}
 \begin{proof}
 We shall show that $U$ is the smallest submodule of $M$ for which $Hom_R(M/N,M/U)=0$. To do this, define $\mathcal{A}=\{B\mid Hom_R(M/N, M/B)=0\}$ and let $U=\cap \mathcal{A}$. Then, $M/U$ can be embedded in $\prod_{B\in \mathcal{A}} M/B$, leading to the conclusion that: 
\begin{align*}
        Hom_R(M/N,M/U)&\leq Hom_R(M/N,\prod_{B\in \mathcal{A}} M/B) \\&\cong \prod_{B\in \mathcal{A}} Hom_R(M/N, M/B) =0. 
    \end{align*}
This shows that $U\in \mathscr{A}$.\
  It is clear that $U\leq B, \forall  B\in \mathscr{A} $.\end{proof}
We will refer to module $_RU$ from Lemma \ref{tot} as the totalizer of $N$ in $M$ and denote it by $Tot^l(N)$.\\

\begin{lemma}\label{BJKNco} The following assertions are equivalent:
    \begin{enumerate}
        \item  For all $L,N $ proper submodules of $_RM$, $L:N$ is a proper submodule of $M$. 
        \item $_RM$ is coprime.
        \item $Hom_R(M/L,M/N) \neq 0$ for proper submodules $L$ and $N$ of $_RM$.
    \end{enumerate}
\end{lemma}
\begin{proof} 

\begin{itemize} 
\item[] 

    \item[$1) \Rightarrow 3) $] Let $L$ and $N$ be  proper submodules of $M$, then $M$ is not equal to $L:N$, meaning that $M/N$ is not equal to $\omega ^{M/L }_0 (M/N)$, which is the intersection of all kernels of homomorphisms from $M/N$ to $M/L$. Therefore, there exists a nonzero homomorphism $f: M/N\rightarrow M/L$. 
    \item[$3)\Rightarrow 1)$]Let \( 0 \neq f \in Hom_R(M/L, M/N) \). If \( f(m + L) \neq 0 + N \), then $m + L \notin \bigcap\{ \ker (g) \mid g: M/L \to M/N \} = \omega_0^{M/N}(M/L) = (N:L)/L$. Consequently, \( (N:L)/L \neq M/L \), which implies \( (N:L) \neq M \) if \( L \) and \( N \) are proper submodules of \( M \).
    \item[$3)\Rightarrow 2)$] If $M/N$ is not equal to zero, denote $T$ as the $\xi(M/N)$-torsion submodule of $M$. If $T$ is not equal to $M$, then $M/T$ would be $\xi(M/N)$-pretorsion-free; however, this would lead to a contradiction since $Hom_R(M/N,M/T)$ would not be equal to zero. Therefore, $T$ must be equal to $M$.
    \item[$2)\Rightarrow 3)$] If $L$ and $N$ are proper submodules of $M$, then $M/N$ is of $\xi (M/L)$-torsion, implying that $Hom_R(M/L,M/N)$ is not equal to zero.
\end{itemize}
\end{proof}

\section{Co-first modules relative to a subclass of $\rpr$.}

 From Proposition 4.3 in \cite{Bic1}, we know that a module $M$ is coprime if, for any two submodules $A$ and $B$ of $M$, the condition $A \Box_M B = M$ implies that either $A = M$ or $B = M$. This is equivalent to saying that $M$ is $\gamma_A^M$-torsion or that $M$ has no nonzero $\gamma_A^M$-torsion quotients.

Similarly, from Proposition \ref{BJKNco}, we have that a module $M$ is coprime if, for any two submodules $A$ and $B$ of $M$, the condition $A : B = M$ implies that either $A = M$ or $B = M$. This is equivalent to stating that $M$ is $\omega_0^{M/A}$-torsion or that $M$ has no nonzero  $\omega_0^{M/A}$-torsion quotients.

We can therefore extend the notion of coprimeness by using a subclass of $\rpr$, as described in the following definition:

\begin{definition}\label{2}
Let $0\neq M\in \rmod$ and $\mathscr{A}\subseteq \rpr$. We say that:
\begin{itemize}
\item[i)] $M$ is  $\mathscr{A}$-co-first if for each $\alpha \in\mathscr{A}$, either $M\in \tfor_\alpha$ or $M$ has no non-zero quotients in $\tfor_\alpha$.
\item[ii)] $M$ is  fully $\mathscr{A}$-co-first if $M$ has no non-zero quotients in $\tfor_\alpha$ for all $ \alpha\in \mathscr{A} $.

\end{itemize}
\end{definition}

    Notice that fully $\mathscr{A}$-co-first modules are $\mathscr{A}$-co-first modules.



\begin{proposition}
    Let $0\neq M\in \rmod$. The following statements are equivalent.
    \begin{itemize}
        \item[(1)] $M$ is $\rrad$-co-first module.
        \item[(2)] $M$ is coprime.
        
    \end{itemize}
\end{proposition}

\begin{proof}
    \begin{itemize}
        \item[] 

        \item[$(1)\Rightarrow (2)$ ] Let $N, K \in \mathscr{L}(M) \setminus \{M\}$. By Remark \ref{tryrej}, $\omega_0^{M/N} \in \rrad$ and $\omega_0^{M/N}(M) \leq N \neq M$, which implies that $M \notin \tfor_{\omega_0^{M/N}}$.
        Since $M$ is $\rrad$-co-first, it has no non-zero quotients in $\tfor_{\omega_0^{M/N}}$. Therefore, for any proper submodule $K$ of $M$, $\omega_0^{M/N}(M/K) \neq M/K$. This guarantees the existence of a nonzero morphism $f: M/K \to M/N$ for each $K, N \lneq M$. Thus, by Lemma \ref{BJKNco}, $M$ is coprime.
        \item[$(2)\Rightarrow (1)$] Let $N\in \mathscr{L}(M)\setminus \{M\}$ and let $\sigma \in \rrad$ be such that $\sigma (M)\neq M$. Then,  by Lemma \ref{BJKNco}, there exists a nonzero morphism $f:M/N\rightarrow M/\sigma(M) $. As  $\sigma\in \rrad$, $\sigma(M/\sigma(M))=0$, and the following diagram is commutative
        \begin{center}
\begin{tikzcd}
 M/N\ar{r}{f}& M/\sigma(M) \\
\sigma(M/N)\ar[hookrightarrow]{u}\ar{r}{f\vert} & \sigma(M/\sigma(M))=0.\ar[hookrightarrow]{u}
\end{tikzcd} 

\end{center}
Thus, $ f(\sigma(M/N))\leq\sigma(M/\sigma(M))=0$, which implies that $\sigma(M/N)\neq M/N $. Therefore, $M$ is $\rrad$-co-first.
\end{itemize}
\end{proof}

\begin{example}\label{7}
 Let $p\in \mathbb{Z}$ be a prime. $_\mathbb{Z}\mathbb{Z}_{p^\infty}$ is $\zpr$-co-first, since for each $N\in\mathscr{L}(\mathbb{Z}_{p^\infty})\setminus \{\mathbb{Z}_{p^\infty}\}$ we have $\mathbb{Z}_{p^\infty}/N\cong \mathbb{Z}_{p^\infty}$. So if $\sigma\in \zpr$ and $_\mathbb{Z}\mathbb{Z}_{p^\infty}\notin \mathbb{T}_\sigma$, then $_\mathbb{Z}\mathbb{Z}_{p^\infty}/N\notin \mathbb{T}_\sigma$ for each proper submodule $N$ of $_\mathbb{Z}\mathbb{Z}_{p^\infty}$, since $\mathbb{T}_\sigma$ is closed under taking isomorphic copies of its elements.
\end{example}

\begin{lemma}
 Every simple module $S\in \rmod$ is $\rpr$-co-first.
\end{lemma}

\begin{proof}
This is clear.
\end{proof}


\begin{lemma}\label{5}
Let $M\in \rmod$ be a non-zero semisimple  module. Then  $M$ is $\rpr$-co-first module  iff  $M\cong S^{(I)}$ for some set $I$ and some simple $R$-module $S$.
\end{lemma}

\begin{proof}
    Let $M\in \rmod$ be a non-zero semisimple module. 
    
    \begin{itemize}
        \item[$\Rightarrow)$] Suppose that $M$ is $\rpr$-co-first module.  If $M\notin \mathbb{T}_{\sigma}$ and $\sigma(M)\neq 0$, then $M=\sigma(M)\oplus N$ with $N\neq 0$, thus $\sigma(N)\leq N\cap \sigma(M)=0$, so $\sigma(N)=0$, thus $\sigma (M)=\sigma(\sigma(M)\oplus N)=\sigma (\sigma(M))\oplus \sigma(N)=\sigma (\sigma(M))$. Hence, $0\neq \sigma(M)\in \mathbb{T}_\sigma $ with  $\sigma(M)\cong M/N$, a contradiction. Therefore, for each $\sigma\in \rpr$, $M\in \mathbb{T}_\sigma$ or $M\in \mathbb{F}_\sigma$.  Since  there exists a simple submodule $S$ of $M$, then $0\neq tr_S(M)$, so $M\in \mathbb{T}_{tr_S}$, therefore  $M\cong S^{(I)}$ for some set $I$.
        \item[$\Leftarrow$)] Suppose that $M\cong S^{(I)}$ for some set $I$ and some simple $R$-module $S$. Then $M\in \mathbb{T}_\sigma $ or $M\in \mathbb{F}_\sigma$. Therefore, if $M\notin \mathbb{T}_\sigma$, then $\sigma(S)=0$, in this case every non-zero quotient of $M$ belongs to $\mathbb{F}_\sigma$.
\end{itemize}
\end{proof}

\begin{definition}
 A submodule $N$ of $M$ is superfluous if $N+K$ is a proper submodule of $M$ for each proper submodule $K$ of $M$, denoted as $N\ll M$.
\end{definition}

Recall Definition 42 of \cite{primepre}.
\begin{definition}
An $R$-module $M$ is called dihollow if every proper fully invariant submodule $N$ of $M$ is superfluous in $M$.
\end{definition}

\begin{proposition}\label{6}
Let $M$ be an $\rpr$-co-first $R$-module. Then $M$ is dihollow.
\end{proposition}

\begin{proof}
Let $_RM$ be a $\rpr$-co-first $R$-module, $N\in  \mathscr{L}_{fi}(M)\setminus \{ M\}$ and $K\in  \mathscr{L}(M)\setminus \{ M\}$. We have $\alpha_N^M(M)=N\neq M$. Since $M$  is $\rpr$-co-first, it follows that $\alpha_N^M(M/K)\neq M/K$.  Moreover, because   $\alpha_N^M\in \rpr$ the following diagram is commutative:
 \begin{center}
\begin{tikzcd}
 M\ar{r}{\pi_K}& M/K \\
\alpha_N^M(M)=N\ar[hookrightarrow]{u}\ar{r}{\pi_K\vert} & \alpha_N^M(M/K).\ar[hookrightarrow]{u}
\end{tikzcd} 

\end{center}
Thus, $ (N+K)/K=\pi_K(N)\leq \alpha_N^M(M/K)\lneq M/K$, which implies that $N+K\lneq M$. Therefore, $N$ is superfluous in $ \mathscr{L}(M)$.
\end{proof}

Proposition \ref{6} shows that $\mathscr{L}_{f.i.}(M)$ is a hollow lattice if $M$ is $\rpr$-co-first. The following example shows that the converse of  Proposition \ref{6} does not hold.

\begin{example}
 $_\mathbb{Z}\mathbb{Z}$ is not $\zpr$-co-first, because  $tr_{\mathbb{Z}_2}(\mathbb{Z}/2\mathbb{Z})=\mathbb{Z}/2\mathbb{Z}$ and $tr_{\mathbb{Z}_2}(\mathbb{Z})=0$. This shows that  $\mathbb{Z}\notin \mathbb{T}_{tr_{\mathbb{Z}_2}}$ although it does have a nonzero quotient in $\mathbb{T}_{tr_{\mathbb{Z}_2}}.$
\end{example}

The following proposition characterizes the $\rpr$-co-first modules.  

\begin{proposition}\label{pf}
The following statements are equivalent:

\begin{itemize}
\item[(1)] $_RM$ is $\rpr$-co-first.
\item[(2)] $M$ is coprime.
\end{itemize}
\end{proposition}
\begin{proof}
\begin{itemize}
\item[] 
    \item[$(1)\Rightarrow (2)$]   Let $N, K\in \mathscr{L}(M)\setminus \{M\}$. Then $\omega_0^{M/N}\in \rpr$ and $\omega_0^{M/N}(M)\leq N\neq M$ (see Remark \ref{tryrej}). Since $M$ is $\rpr$-co-first, $M\notin \mathbb{T}_{\omega_0^{M/N}}$ and $ K\lneq M$, then $\omega_0^{M/N}(M/K)\neq M/K$. Thus, there exists a nonzero morphism $f:M/K\rightarrow M/N$. Therefore, $M$ is coprime according to Lemma \ref{BJKNco}.
    \item[$(2)\Rightarrow (1)$]  Let $N\in \mathscr{L}(M)\setminus \{M\}$ and $\sigma \in \rpr$ such that $\sigma(M/N)=M/N$. As $M$ is  coprime, by Lemma \ref{BJKNco}, there exists a set $X$ and an epimorphisms $g:(M/N)^{(X)}\rightarrow M$.  Thus:\begin{center}
$M=g((M/N)^{(X)})=g((\sigma(M/N))^{(X)})=g(\sigma((M/N)^{(X)}))\leq \sigma(M)\leq M$
        \end{center}
        Therefore, $\sigma(M)=M$, so $M$ is $\rpr$-co-first.
        \end{itemize}
\end{proof}

\begin{example}\label{ejemcovtcf}
    Let $p \in \mathbb{Z}$ be a prime. Consider $M = \mathbb{Z}_{p^2}$. Let $K \lneq M$ and $I$ be an ideal of $\mathbb{Z}$.
    Then $M/K \cong \mathbb{Z}_{p^i}$ with $i \in \{1,2\}$ and there exists $n \in \mathbb{Z}$ such that $I = n\mathbb{Z}$.
    If $I(M/N) = M/N$, then $I + p^i \mathbb{Z} = \mathbb{Z}$, therefore $n$ and $p^i$ are coprime, so $n$ and $p$ are coprime, $I + p \mathbb{Z} = \mathbb{Z}$, then $IM = M$. Thus $M$ is an $\trad$-co-first module.\\
    On the other hand, $\mathbb{Z}_p$ is a nonzero quotient  of $\mathbb{Z}_{p^2}$ that does not generate a $\xi(\mathbb{Z}_{p^2})$, so $\mathbb{Z}_{p^2}$ is not a coprime module.
    \end{example}
\section{The classes of $\mathscr{A}$-co-first modules and of fully $\mathscr{A}$-co first modules.}

In this section, we study the classes of $\sigma$-fully co-first modules and examine the closure properties they satisfy. We will begin this section by introducing the lattice of conatural and cohereditary classes.

Recall that $R$-quot is the class of classes of modules  which are closed under quotients.

\begin{proposition}[\cite{OTLONACC}, Proposition 16]\label{seudocompinrquot}
    If $\mathscr{Q}\in R$-quot, then $\mathscr{Q}$ has an unique pseudocomplement $\mathscr{Q}^{\perp_\twoheadrightarrow}$ within $R$-quot, which can be described as the  class $\mathscr{K}$ comprising all modules $N$ with no non-zero quotients in $\mathscr{Q}$.
  
\end{proposition}

Recall that the skeleton of a lattice $\mathscr{L}$ is the set of the pseudocomplements of its elements. The skeleton of $R$-quot is called $\rconat$ and its elements are called conatural classes. These are characterized in Theorem \ref{CNCNC} below.

\begin{definition}\label{CNC}[\cite{OTLONACC}, Definition 22]
  Let $\mathscr{C}$ be a class of $R$-modules. We say that $\mathscr{C}$ satisfies condition $(CN)$ if
 \begin{center}
 $(\forall 0\neq f:M \twoheadrightarrow N \text{ epimorphism, such that there exist epimorphisms }  N \twoheadrightarrow K \twoheadleftarrow C\textit{ with } C \in \mathscr{C} \text{ and }   0 \neq K ), M \in \mathscr{C}.$
\end{center}

(If every non-zero quotient of $M$ shares a non-zero quotient with some  $C\in \mathscr{C}$, then $M\in \mathscr{C})$.
\end{definition}

\begin{theorem}[\cite{OTLONACC}, Theorem 23]\label{CNCNC}
    Let $\mathscr{C}$ be a class of $R$-modules. The following conditions are equivalent:
    \begin{itemize}
        \item[(1)] $\mathscr{C}\in \rconat$,
        \item[(2)] $\mathscr{C}$ satisfies condition $(CN)$,
        \item[(3)] $\mathscr{C} \in R$-quot and $(\mathscr{C}^{\perp_\twoheadrightarrow})^{\perp_\twoheadrightarrow} = \mathscr{C}$.
        \item[(4)] $\mathscr{C}=\mathscr{D}^{\perp_\twoheadrightarrow}$ for some $\mathscr{D}\in R$-quot.
    \end{itemize}
\end{theorem}

We will use the following facts about conatural classes:
\begin{remark}[\cite{OTLONACC}, Propositions 26 and 28]\label{remarkconat}
    Let $\mathscr{C}$ be a subclass of $\rmod$. The following statements hold:
    \begin{enumerate}
        \item $\{0\}$ and $\rmod$ are conatural classes.
        \item If $\mathscr{C}\in R$-quot, then:\begin{center}

$(\mathscr{C}^{\perp_\twoheadrightarrow})^{\perp_\twoheadrightarrow} = \{ M \in \rmod \mid \text{ for every non-zero epimorphism } f: M \rightarrow N, \text{ there exists an epimorphism } g: N \rightarrow K \text{ with } 0\neq K \in \mathscr{C} \}$            
        \end{center}
is the smallest conatural class containing $\mathscr{C}$. In particular, if $\mathscr{C} \in \rconat$, then $(\mathscr{C}^{\perp_\twoheadrightarrow})^{\perp_\twoheadrightarrow} = \mathscr{C}$.
        \item $\rconat$ is a Boolean lattice.
    \end{enumerate}
    \end{remark}

\begin{definition}\label{defclasses}
 For a subset $\mathscr{A} \subseteq \rpr$, we will use the following classes: \begin{itemize} 
 \item $\mathbb{T}_\mathscr{A}=\{ M\in \rmod\vert M\in \mathbb{T}_\sigma\textit{, } \forall\sigma\in \mathscr{A}\}$.
  \item $\mathbb{F}_\mathscr{A}=\{ M\in \rmod\vert M\in \mathbb{F}_\sigma\textit{, } \forall\sigma\in \mathscr{A}\}$.
 \item $\mathscr{P}_{\mathscr{A}}$,  the class of all fully $\mathscr{A}$-co-first modules
.
\item $\overline{\mathbb{P}}_{\mathscr{A}}$,  the class that includes all $\mathscr{A}$-co-first modules and the zero module. \end{itemize} We will use $\mathscr{P}_{\sigma}$ and $\overline{\mathbb{P}}_{\sigma}$ as shorthand for $\mathscr{P}_{\{\sigma\}}$ and $\overline{\mathbb{P}}_{\{\sigma\}}$, respectively.
\end{definition}

\begin{remark}\label{cop8}
For any $\sigma,\beta\in \rpr$ and $\mathscr{A}\subseteq \rpr$, the following statements hold:
\begin{itemize}
\item[(1)]  $\mathscr{P}_{\mathscr{A}}=\bigcap\limits_{r\in \mathscr{A}}\mathscr{P}_{r}$ and 
 $\overline{\mathbb{P}}_{\mathscr{A}}=\bigcap\limits_{r\in \mathscr{A}}\overline{\mathbb{P}}_{r}$.
 \item[(2)] $\overline{\mathbb{P}}_{\sigma}=\mathscr{P}_{\sigma} \cup \mathbb{T}_{\sigma}$.
 \item[(3)]  $\mathscr{P}_{\mathscr{A}}\subseteq \overline{\mathbb{P}}_{\mathscr{A}}$ and 
 $\mathbb{T}_{\mathscr{A}}\subseteq \overline{\mathbb{P}}_{\mathscr{A}}$.
 \item[(4)] If $\sigma \preceq \beta$, then $\mathscr{P}_\beta\subseteq \mathscr{P}_\sigma$.
\end{itemize}
\end{remark}

\begin{proposition}\label{propPsig}

    Let $\sigma \in \rpr$. Then $\mathscr{P}_\sigma=(\mathbb{T}_\sigma)^{\perp_\twoheadrightarrow}$. 
    
\end{proposition}

\begin{proof}
  
Let $\sigma \in \rpr$. Then, by Proposition \ref{seudocompinrquot}, we have
\begin{align*}
    \mathscr{P}_\sigma &= \{ M \in \rmod \mid M \text{ has no non-zero quotients in } \tfor_\alpha \} \\
    &= \{ M \in \rmod \mid \text{ the only } \sigma\text{-pretorsion quotient of } M \text{ is } 0 \} \\
    &= (\mathbb{T}_\sigma)^{\perp_{\twoheadrightarrow}}.
\end{align*}

\end{proof}

\begin{remark}
    Let $\sigma \in \rpr$. Note that since $\mathbb{T}_\sigma=\mathbb{T}_{\widehat{\sigma}}$,  it follows from Proposition \ref{propPsig} that $\mathscr{P}_\sigma=\mathscr{P}_{\widehat{\sigma}}$.
\end{remark}

\begin{corollary}\label{conatp}
    Let $\sigma \in \rpr$. Then $\mathscr{P}_\sigma$ is a conatural class. 
\end{corollary}
\begin{proof}
    Let $\sigma \in \rpr$. We have that $\mathbb{T}_\sigma\in R$-quot, and By Proposition \ref{propPsig} $\mathscr{P}_\sigma=(\mathbb{T}_\sigma)^{\perp_\twoheadrightarrow}$. We conclude using Theorem \ref{CNCNC}, (4).
\end{proof}

\begin{proposition}
    Let $N$ be a fully invariant proper submodule of the module $_RM$. If $M\in \mathscr{P}_{\alpha_N^M}$, then $N$ is superfluous in $M$. 
\end{proposition}

\begin{proof}
Let $N\in \mathscr{L}_{fi}(M)\setminus \{ M\}$. If $M\in \mathscr{P}_{\alpha_N^M}$, then for all  $ K\lneq M$, we have $\alpha_N^M(M)=N\neq M$, and  consequently $\alpha_N^M(M/K)\neq M/K$.
Let us denote $\pi:M\to M/K$ the canonical projection.  Then $(N+K)/K=\pi(N)=\pi(\alpha_N^M(M))\leq \alpha_N^M(M/K)\lneq M/K$. Hence $N+K\lneq M$.  Therefore, $N$ is superfluous in $M$. 
\end{proof}

\begin{remark}
    We have that $\mathscr{P}_{(\_)}:\rpr \to \rconat$ is an antimorphims of posets by Remark \ref{cop8} and  Corollary \ref{conatp}.
\end{remark}

\begin{corollary}\label{24}
The following statements are  equivalent for $\sigma  \in R$-pr:
\begin{itemize}
    \item[(1)] $\mathscr{P}_\sigma=\rmod$.
    \item[(2)] $\widehat{\sigma}= \underline{0}$.
    \item[(3)] $\mathbb{T}_{\sigma}=\{0\}$.
\end{itemize}
\end{corollary}

\begin{proof}
\begin{itemize}
\item[]
   \item[$(1)\Rightarrow (3)$]
Since $\rmod=\mathscr{P}_\sigma=(\mathbb{T}_{\sigma})^{\perp_\twoheadrightarrow}$, then  $\mathbb{T}_\sigma\subseteq ((\mathbb{T}_\sigma)^{\perp_\twoheadrightarrow})^{\perp_\twoheadrightarrow}=(\mathscr{P}_\sigma)^{\perp_\twoheadrightarrow}=\rmod^{\perp_\twoheadrightarrow}=\{0\}$.
   \item[$(3)\Rightarrow (1)$] Since $\mathbb{T}_\sigma=\{0\}$, then $\mathscr{P}_\sigma =(\mathbb{T}_\sigma) ^{\perp_\twoheadrightarrow}=\{0\}^{\perp_\twoheadrightarrow}=\rmod$.
    \item[$(2)\Rightarrow (3)$]  This is because $\mathbb{T}_\sigma=\mathbb{T}_{\widehat{\sigma}}$.
    \item[$(3)\Rightarrow (2)$] Since $\widehat{\sigma}\in\rid$, $\widehat{\sigma}(M)\in \mathbb{T}_{\widehat{\sigma}}=\mathbb{T}_\sigma=\{0\}$, for each $M\in \rmod$. \end{itemize}
   
\end{proof}


\begin{definition}
    Let $\mathscr{A}$ be a subclass of $\rpr$. Let us define
    \begin{center}
        $\mathcal{CN}_\mathscr{A}=\{\mathscr{P}_\sigma\vert \sigma \in \mathscr{A}\}$.
    \end{center}
\end{definition}

\begin{proposition}
    Let us denote  $\rpr^*=\rpr\setminus \{ \underline{0}\}$. The following statements are equivalent:
    \begin{itemize}
    \item[(1)] $\rmod\notin\mathcal{CN}_\mathscr{\rpr^*}$.
    \item[(2)]  $\mathbb{T}_\sigma\neq \{ 0\}$ for each  $\sigma \in \rpr^*$.
    \item[(3)] $R$ is a left $V$-ring.
\end{itemize}
\end{proposition}

\begin{proof}
\begin{itemize}
\item[] 

\item[(1)$\Leftrightarrow$(2)] It  follows from the Corollary \ref{24}.
    \item[(3)$\Rightarrow$(2)] Let  $R$ be a $V$-ring and $\sigma\in R-pr\setminus \{ \underline{0}\}$. By Lemma 6 of \cite{pre}, there exists  $S\in \rsimp$ such that  $\sigma(E(S))\neq 0$, so $\alpha_S^{E(S)}\preceq \alpha_{\sigma(E(S))}^{E(S)}\preceq \sigma$. As  $R$ is a  $V$-ring, $S$ is an injective module, thus   $\alpha_S^{E(S)}=\alpha_S^S$, and then  $\alpha _S^S(S)=S\leq \sigma(S)\leq S$. Hence $S\in \mathbb{T}_\sigma$, so $\mathbb{T}_\sigma\neq \{ 0\}$.
    \item[(2)$\Rightarrow$(3)] Assume $\mathbb{T}_\sigma\neq \{ 0\}$ for all $\sigma \in R\text{-pr}\setminus \{ \underline{0}\}$. Given  $S\in R-$simp, $\alpha_S^{E(S)}\neq \underline{0}$, then there exists a nonzero $M\in\rmod$ such that  $\alpha_S^{E(S)}(M)=M$.  Thus  $0\neq
\alpha _M^M\preceq \alpha_S^{E(S)}$. Since  $\alpha_S^{E(S)}$ is an  atom in $\rpr$, by Theorem 7 of \cite{pre}, we have $\alpha _M^M=\alpha_S^{E(S)}$.\\ Therefore, $ \alpha_S^{E(S)}$ is an idempotent preradical for each  $S\in \rsimp$, making $R$ a  $V$-ring by Theorem 17 of \cite{pre}.
\end{itemize}

\end{proof}

\begin{proposition}\label{8.3}
Let $\sigma \in \trad$. Then $\mathscr{P}_\sigma=\{M\in \rmod \mid \sigma(M)\ll M\}.$
\end{proposition}
 \begin{proof}
Let $M\in \rmod$ and $N\leq M$. Then $$\sigma (\frac{M}{N})=\sigma(R)\cdot\frac{M}{N}=\frac{\sigma(R)\cdot(M)+N}{N}=\frac{\sigma(M)+N}{N},$$ because $\sigma$ is a $t$-radical (see \cite{Es} Exercise VI,5). Now,  if $M\in \mathscr{P}_\sigma$, then its only $\sigma$-pretorsion quotient is $0$. $\sigma(M)+ N=M$ implies   $$\frac{M}{N}=\frac{\sigma(M)+ N}{N}=\sigma (\frac{M}{N}).$$ Hence $N=M$, so $\sigma(M)$ is superfluous in $M$.\\
Conversely, if $\sigma(M)$ is superfluous in $ M$, then  $\sigma(M)+ N=M$, implies  $N=M$. Thus the only $\sigma$-pretorsion quotient of $M$ is $0$ and therefore $M\in \mathscr{P}_\sigma$.
 \end{proof}

The following proposition found in \cite{OTLONACC} contains a demonstration to make it easy to read.
\begin{proposition}\label{cftrad}
    Let $\sigma \in \trad$. Then $\mathbb{T}_\sigma\in \rconat$
\end{proposition}

\begin{proof}
\begin{itemize}
    \item[]
    \item[$\subseteq )$] Let $M\in (\mathscr{P}_\sigma)^{\perp_\twoheadrightarrow} $. Then $M/\sigma(M)\in (\mathscr{P}_\sigma)^{\perp_\twoheadrightarrow}$. If $\sigma (M)\lneq M$, then $M/\sigma(M)\neq 0$.  Since $\sigma$ is a radical, then $\sigma(M/\sigma(M))=0$, which is a superfluous submodule of $M/\sigma(M) $. So $0\neq M/\sigma(M)\in \mathscr{P}_\sigma $, by Proposition \ref{8.3}.  Hence  $0\neq M/\sigma(M)\in \mathscr{P}_\sigma \cap(\mathscr{P}_\sigma)^{\perp_\twoheadrightarrow}$ which is a contradiction. Therefore $(\mathscr{P}_\sigma)^{\perp_\twoheadrightarrow}\subseteq \mathbb{T}_\sigma$.
    \item[$\supseteq )$] Suppose that $M\notin (\mathscr{P}_\sigma)^{\perp_\twoheadrightarrow}$. Then there exists $N\lneq M$ such that $M/N\in \mathscr{P}_\sigma$, so $\sigma(M/N)\neq M/N$. Hence $\sigma(M)\lneq M$, i.e.  $M\notin \mathbb{T}_\sigma$. Therefore $\mathbb{T}_\sigma\subseteq (\mathscr{P}_\sigma)^{\perp_\twoheadrightarrow}$.
\end{itemize}

\end{proof}

\begin{lemma}\label{coroL}
    Let $\sigma\in \trad$. Then $\mathbb{F}_\sigma\subseteq \mathscr{P}_\sigma$.
\end{lemma}

\begin{proof}
  By \cite{OTLONACC}, Proposition 28,  $\rconat$ is a Boolean lattice. By Proposition \ref{cftrad}, $\mathbb{T}_\sigma=(\mathscr{P}_\sigma)^{\perp_\twoheadrightarrow}$ is a conatural class, so $\mathscr{P}_\sigma=(\mathbb{T}_\sigma)^{\perp_\twoheadrightarrow}$. As $\sigma\in \trad$, $\mathbb{F}_\sigma$ is closed under quotients. Clearly $\mathbb{F}_\sigma\cap\mathbb{T}_\sigma=\{0\}$, hence $\mathbb{F}_\sigma\subseteq (\mathbb{T}_\sigma)^{\perp_\twoheadrightarrow}=\mathscr{P}_\sigma$.
  
\end{proof}

Later, in Example \ref{ejem4.17}, we will observe that the inclusion can be proper.

\begin{proposition}
    Let $\sigma\in\rrid$. Then $\sigma\in \trad$ if and only if $\mathbb{F}_\sigma\subseteq \mathscr{P}_\sigma$.
\end{proposition}

\begin{proof}
\begin{itemize}
    \item[]
    \item[$\Rightarrow$)] If $\sigma\in \trad$, according to Lemma \ref{coroL}, we have $\mathbb{F}_\sigma\subseteq \mathscr{P}_\sigma$. \item[$\Leftarrow$)] Now, assume $\mathbb{F}_\sigma\subseteq \mathscr{P}_\sigma$ and consider an epimorphism $g:M\to N$ where $M$ is an element of $\mathbb{F}_\sigma$. Let $K=Ker(g)$. Given that $\sigma$ is idempotent, $\sigma(M/K)=H/K$ is a $\sigma$ pretorsion module, hence $\sigma(H/K)=H/K$. This means $H/K$ belongs to $\mathbb{T}_\sigma$. Furthermore, since $H\leq M\in \mathbb{F}_\sigma$, it follows that $H\in \mathbb{F}_\sigma$, leading to $H\in \mathscr{P}_\sigma=(\mathbb{T}_\sigma)^{\perp_\twoheadrightarrow}$, as established in Proposition \ref{propPsig}. Consequently, $H$ has no nonzero quotients in $\mathbb{T}_\sigma$, which implies $H/K=0$. Therefore, $\sigma(N)\cong\sigma(M/K)=H/K=0$, indicating that $N\in \mathbb{F}_\sigma$. Thus, $\mathbb{F}_\sigma$ is closed under quotients, it follows that $\sigma \in \trad$, according to Proposition I.2.8 from \cite{Bic}.
\end{itemize}

\end{proof}

\begin{proposition}\label{propseudo}
    Let $\mathscr{C}\in\rconat$. If there exists an idempotent $t$-radical $\sigma$ such that $\mathscr{C}^{\perp_\twoheadrightarrow}=\mathscr{P}_{\sigma}$ then  $\mathscr{C}$ is closed under coproducts.
\end{proposition}

\begin{proof}
     Let $\sigma$ be an idempotent $t$-radical such that $\mathscr{C}^{\perp_\twoheadrightarrow}=\mathscr{P}_{\sigma}$. Since $\rconat$ is a Boolean lattice then $\mathscr{C}=(\mathscr{C}^{\perp_\twoheadrightarrow})^{\perp_\twoheadrightarrow}$, and by Proposition \ref{cftrad} $\mathscr{P}_{\sigma}^{\perp_\twoheadrightarrow}=\mathbb{T}_\sigma$.
     Thus $\mathscr{C}=\mathbb{T}_\sigma$, therefore, $\mathscr{C}$ is closed under coproducts.
\end{proof}

%
\begin{theorem}[Theorem 30, \cite{OSLOMC}]\label{ECNCUC}
    $\rconat\subseteq \rTORS$  $\Longleftrightarrow R \text{ is a left MAX ring} $.
\end{theorem}

\begin{theorem}\label{MAXRing}
   If $\mathcal{CN}_\trad=\rconat$
      then  $R$ is a left MAX ring.
\end{theorem}\begin{proof}
   Let $\zeta\in\rconat$. Then $\zeta^{\perp_\twoheadrightarrow} \in \rconat$, so there exists 
      $\sigma \in \trad$ such that 
        $\zeta^{\perp_\twoheadrightarrow} = \mathscr{P}_{\sigma}$. 
     Then $\zeta$ is closed by coproducts by Proposition \ref{propseudo}. Therefore $R$ is a left MAX ring by Theorem \ref{ECNCUC}.
     
\end{proof}

\begin{lemma}\label{tpcfq}
    Let $R$ be a perfect ring and $\sigma \in \rrid$. $\mathbb{T}_\sigma$ is closed under projective covers if and only if $\mathbb{F}_\sigma$ is closed under quotients.

\end{lemma}

\begin{proof}

    \begin{itemize}

    \item[]

    \item[$\Rightarrow ) $] Let $g:M\to N$ be an epimorphism with $M\in \mathbb{F}_\sigma$ and $K=g^{-1}(\sigma(N))$. Since $R$ is a perfect ring, there exists a projective cover $P\twoheadrightarrow \sigma(N)$, so there exists a morphism $h:P\to K$ which makes the following diagram commutative:
    \begin{center}
        \begin{tikzcd}
            & P\ar[d,twoheadrightarrow]\ar[dl,"h"]\\
            K\ar[r,twoheadrightarrow, "g\vert"'] & \sigma(N).
        \end{tikzcd}
    \end{center}
    Now, since $\sigma$ is idempotent, $\sigma(N)\in \mathbb{T}_\sigma$ and, since $\mathbb{T}_\sigma$ is closed under projective covers, $P\in \mathbb{T}_\sigma$. Then $h(P)\in \tfor_\sigma$. On the other hand, $h(P)\leq K\leq M\in \mathbb{F}_\sigma$, so $h(P)\in \mathbb{F}_\sigma$. Thus, $h(P)=0$.\\
    Thus $\sigma(N)=gh(P)=0$, which implies that $N\in \mathbb{F}_\sigma$.
    
     \item[$\Leftarrow)$] Let $M \in \mathbb{T}_\sigma$, and 
\begin{tikzcd} 
P \ar[twoheadrightarrow]{r}{\pi} & M 
\end{tikzcd}
a projective cover. Consider the following diagram:
\begin{center}
    \begin{tikzcd}
        \sigma(P) \ar[hook]{d} \ar[twoheadrightarrow]{r}{\pi_{\mid}} & \pi\sigma(P) \ar[hook]{d} \\
        P \ar[twoheadrightarrow]{d}{p_1} \ar[twoheadrightarrow]{r}{\pi} & M \ar[twoheadrightarrow]{d}{p_2} \\
        P/\sigma(P) \ar[twoheadrightarrow]{r}{\overline{\pi}} & M/\pi\sigma(P),
    \end{tikzcd}
\end{center}
where $P/\sigma(P) \in \mathbb{F}_{\sigma}$, because $\sigma$ is a radical. Then, by hypothesis, $M/\pi\sigma(P)$ also belongs to $\mathbb{F}_{\sigma}$. Notice that it is also in $\mathbb{T}_{\sigma}$, because it is a quotient of $M$. Therefore, $M/\pi\sigma(P) =_R 0$. Thus, $P = \sigma(P)$. Hence, $P \in \mathbb{T}_\sigma$.
\end{itemize}

\end{proof}

\begin{theorem}\label{conatperfect1}
    If $R$ is a left perfect ring, then $\rconat= \{\tfor_\sigma \hspace{2mm}\vert\hspace{2mm} \sigma\in \trad\}$. 
\end{theorem}
\begin{proof}
Suppose $R$ is a left perfect ring.

$\subseteq)$ Let  $\mathscr{C}$ be a conatural class. As a left perfect ring is a left MAX ring, then  by Theorem \ref{ECNCUC},  each conatural class is the pretorsion class of an idempotent radical $\sigma$, hence  $\mathscr{C}=\tfor_\sigma$. Now we will see that $\sigma\in \trad$. As each conatural class is closed under projective covers and $R$ is a left perfect ring, then  $\mathbb{F}_\sigma$ is closed under quotients by Lemma \ref{tpcfq}. Then $\sigma\in\trad$.

$\supseteq)$ If $\sigma\in\trad$ then $\mathbb{T_\sigma}$ is a conatural class, by Proposition \ref{cftrad}.

\end{proof}

\begin{corollary}\label{conatforleftperfect}
    If $R$ is a left perfect ring, then $\mathcal{CN}_\trad=\rconat$.
\end{corollary}

\begin{proof}
    Since $R$ is a perfect ring,  $\rconat= \{\tfor_\sigma \hspace{2mm}\vert\hspace{2mm} \sigma\in \trad\}$ according to Theorem \ref{conatperfect1}. As $\rconat$ is a Boolean lattice by Remark \ref{remarkconat}, we have that $\rconat =\{\zeta^{\perp_\twoheadrightarrow}\vert \zeta\in\rconat\}$. Therefore
    \begin{align*}
        \rconat&=\{ \{(\tfor_\sigma)^{\perp_\twoheadrightarrow} \hspace{2mm}\vert\hspace{2mm} \sigma\in \trad\}\\
        &= \{\mathscr{P}_\sigma \hspace{2mm}\vert\hspace{2mm} \sigma\in \trad\}
        \\&=\mathcal{CN}_\trad.
    \end{align*}
\end{proof}

\begin{corollary}
   
For a semilocal ring \( R \), the following statements are equivalent:
\begin{enumerate}
    \item \( R \) is a left perfect ring.
    \item \( \mathcal{CN}_\trad = \rconat \).
\end{enumerate}

\end{corollary}
\begin{proof}
\begin{itemize}
    \item[]
    \item[$1)\Rightarrow 2)$] This follows from Corollary \ref{conatforleftperfect}.
\item[$2)\Rightarrow 1)$ ] If \( \mathcal{CN}_\trad = \rconat \), then $R$ is a left MAX ring by Theorem \ref{MAXRing}. Since $R$ is a semilocal ring, then $R$ is a left perfect ring, by Bass Theorem P (\cite{An}, Theorem 28.4).
    
\end{itemize} 
    
\end{proof}



\section{$\mathscr{A}$-second modules vs $\mathscr{A}$-co-first modules}

In this final section, we compare the extension of the notion of coprimeness proposed in this article with the extension of the notion of secondness introduced in \cite{Sec}, as well as the conditions under which these two notions coincide.

\begin{definition}[\cite{Sec} Definition 7.1]
Let $M$ be a non-zero $R$-module and $\mathscr{A} \subseteq \rpr$. We say that $M$ is $\mathscr{A}$-second if for every $\alpha \in \mathscr{A}$:
\begin{center}
    $\alpha(M) = 0$ or $\alpha(M) = M$.
\end{center}
\end{definition}

For $\mathscr{A} \subseteq \rpr$, we will denote by $\mathbb{S}_{\mathscr{A}}$ the class of all $\mathscr{A}$-second modules and the zero module.
In particular, for $\sigma \in \rpr$, we will use $\mathbb{S}_{\sigma}$ instead of $\mathbb{S}_{\{\sigma\}}$.

\begin{remark}[See \cite{Sec}, Definition 7.2, and recall Definition \ref{defclasses}]\label{8s}
For any $\sigma \in \rpr$ and $\mathscr{A} \subseteq \rpr$, we have that: 
\begin{enumerate}
    \item $\mathbb{S}_{\mathscr{A}} = \bigcap\limits_{r \in \mathscr{A}} \mathbb{S}_{r}$.
    \item $\mathbb{S}_{\sigma} = \mathbb{T}_{\sigma} \cup \mathbb{F}_{\sigma}$.
    \item $\mathbb{T}_{\mathscr{A}} \subseteq \mathbb{S}_{\mathscr{A}}$.
    \item $\mathbb{F}_{\mathscr{A}} \subseteq \mathbb{S}_{\mathscr{A}}$.
\end{enumerate}  
\end{remark}

\begin{theorem}\label{tsvc1}
    Let $\mathscr{A} \subseteq \trad$. Every $\mathscr{A}$-second module is an $\mathscr{A}$-co-first module.
\end{theorem}

\begin{proof}
By Corollary \ref{coroL}, we have  $\mathbb{F}_\sigma\subseteq \mathscr{P}_\sigma$. Furthermore, according  to Remark \ref{8s}, (2),  $\mathbb{S}_\sigma=\mathbb{T}_\sigma\cup\mathbb{F}_\sigma$. Besides, Remark \ref{cop8}, (2), $\overline{\mathbb{P}_\sigma}=\mathbb{T}_\sigma\cup \mathscr{P}_\sigma$. Consequently,  $\mathbb{S}_\sigma\subseteq \overline{\mathbb{P}_\sigma}$. Therefore, $\mathbb{S}_\mathscr{A}\subseteq  \overline{\mathbb{P}_\mathscr{A}}$ by Remark \ref{8s}, (1) and Remark \ref{cop8}, (1).
\end{proof}

The converse of Theorem \ref{tsvc1} generally does not hold, as the following counterexample shows.

\begin{example}\label{ejemsvc}
    Let $p \in \mathbb{Z}$ be a prime. Consider $M = \mathbb{Z}_{p^2}$. By Example \ref{ejemcovtcf}, $M$ is an $\mathbb{Z}\textbf{-trad}$-co-first module.\\
    On the other hand, $p\mathbb{Z}$ is an ideal of $\mathbb{Z}$ and $p\mathbb{Z}M \cong \mathbb{Z}_p$, so $M$ is not an $\mathbb{Z}\textbf{-trad}$-second module. Note  that a $\mathbb{Z}\textbf{-trad}$ second module is the same as a second module in the usual sense.
\end{example}

\begin{example}\label{ejem4.17}
In general, the converse of Lemma  \ref{coroL} is not true. Example \ref{ejemsvc} shows that if $\sigma = p\mathbb{Z} \cdot \_$ and $M = \mathbb{Z}_{p^2}$, then $M \in \overline{\mathbb{P}}_\sigma$, $M \notin \mathbb{L}_\sigma$.
\end{example}

\begin{proposition}\label{psvc1}
    Let $\sigma \in \rrad$. If $\overline{\mathbb{P}}_\sigma = \mathbb{S}_\sigma$, then $\sigma \in \trad$.
\end{proposition}

\begin{proof}
    Let $\sigma \in \rrad$ such that $\overline{\mathbb{P}}_\sigma = \mathbb{S}_\sigma$. Then $\mathscr{P}_\sigma = \mathbb{F}_\sigma$ by Remark \ref{8s} and Remark \ref{cop8}. Now, by Proposition \ref{propPsig}, $\mathscr{P}_\sigma = (\mathbb{T}_\sigma)^{\perp_\twoheadrightarrow}$, therefore $\mathbb{F}_\sigma = (\mathbb{T}_\sigma)^{\perp_\twoheadrightarrow}$. Thus, $\mathbb{F}_\sigma$ is closed under quotients, so $\sigma \in \trad$.
\end{proof}

\begin{proposition}\label{radsec=radcop}
   Let $R$ be a ring. The following statements are equivalent:
   \begin{itemize}
       \item[(1)] $\overline{\mathbb{P}}_{rad} = \mathbb{S}_{rad}$.
   \item[(2)] $R$ is a left $V$-ring.
   \end{itemize}
\end{proposition}

\begin{proof}
    \begin{itemize}
    \item[]
       \item[ (1)$\Rightarrow$ (2)] By Proposition \ref{psvc1}, \( rad \in \trad \). Remark \ref{8s} and Remark \ref{cop8} imply \( \mathscr{P}_{rad} = \mathbb{F}_{rad} \), so \( \mathbb{F}_{rad} \in \rconat \). For each left ideal \( I \lneq R \), there exists a maximal left ideal \( \mathscr{M} \) of \( R \) such that \( I \leq \mathscr{M} \). Thus, each nonzero quotient of \( R \) has a simple quotient, and since every simple module \( S \)  is \( rad \)-torsion free, we have \( R \in (\mathbb{F}_{rad}^{\perp_\twoheadrightarrow})^{\perp_\twoheadrightarrow} \). By Remark \ref{remarkconat}, \( R \in \mathbb{F}_{rad} = (\mathbb{F}_{rad}^{\perp_\twoheadrightarrow})^{\perp_\twoheadrightarrow} \), leading to \( rad(R) = 0 \). Furthermore, since \( rad \in \trad \) by Proposition \ref{psvc1}, it follows that \( rad(M) = rad(R)M = 0 \). Therefore, \( R \) is a left \( V \)-ring.

\item[ (1) $\Rightarrow$ (2)] Since \( R \) is a left \( V \)-ring, we have \( rad(M) = 0 \)  for every \( M \in \rmod \), by exercise 13.10 of \cite{An}. Thus \( \mathbb{F}_{rad} = \rmod \) and \( \mathbb{T}_{rad} = \{0\} \). Consequently, \( \mathbb{T}_{rad}^{\perp_\twoheadrightarrow} = \rmod \). Therefore, \( \overline{\mathbb{P}}_{rad} = \mathbb{S}_{rad} \) by Proposition \ref{propPsig}, Remark \ref{8s}, and Remark \ref{cop8}.
    \end{itemize}
\end{proof}

\begin{proposition}
   Let $R$ be a ring. The following statements are equivalent:
    \begin{itemize}
       \item[(1)] $\overline{\mathbb{P}}_\sigma = \mathbb{S}_\sigma$ for every $\sigma \in \rpr$.
   \item[(2)] $R$ is a semisimple ring.
    \end{itemize}
\end{proposition}

\begin{proof}
   \begin{itemize}
    \item[] 
     
\item[(1)$\Rightarrow$(2)] Let \( S \in \rsimp \). Then there exists a maximal left ideal \( I \) of \( R \) such that \( R/I \cong S \). This implies that \( R \notin \mathbb{T}_{\alpha_S^S}^{\perp_\twoheadrightarrow} = \mathscr{P}_{\alpha_S^S} \). Since \( \overline{\mathbb{P}}_{\alpha_S^S} = \mathbb{S}_{\alpha_S^S} \), it follows that \( \mathscr{P}_{\alpha_S^S} = \mathbb{F}_{\alpha_S^S} \), as noted in Remark \ref{8s} and Remark \ref{cop8}. Consequently, \( R \notin \mathbb{F}_{\alpha_S^S} \). Since \( 0 \neq \alpha_S^S(R) \), it follows that \( S \) embeds into \( R \). Therefore, \( R \) is a left Kasch ring. Furthermore, by Proposition \ref{radsec=radcop}, \( R \) is a left \( V \)-ring, which implies that each \( S \in \rsimp \) is an injective module and a direct summand of \( R \). Consequently, each \( S \in \rmod \) is a projective module.

Now, suppose \( \text{soc}(R) \neq R \). Then there exists a maximal left ideal \( I \) of \( R \) such that \( \text{soc}(R) \leq I \). Thus \( R/I \) is a projective simple module, which means that the following exact sequence splits:

\begin{center}
    \begin{tikzcd}
        0 \ar[r] & I \ar[hookrightarrow, r] & R \ar[twoheadrightarrow, r] & R/I \ar[r] & 0.
    \end{tikzcd}
\end{center}

Therefore, there exists a simple submodule \( S \) of \( R \) such that \( R = I \oplus S \). In particular, \( 0 = I \cap S  \geq \text{soc}(R) \cap S=S \), which is a contradiction. Hence, we conclude that \( \text{soc}(R) = R \). This means that \( R \) is a semisimple ring.

        \item[(2)$\Rightarrow$(1)]  Let $\sigma\in \rpr$. By Remarks \ref{cop8} and \ref{8s}, it is sufficient to prove that $\mathscr{P}_\sigma=\mathbb{F}_\sigma$. Let $0\neq M \in \rmod$. Since $R$ is a semisimple ring, there exists a set $I$ and a family $\{S_i\}_{i\in I}$ of simple modules such that  $M\cong \bigoplus\limits _{i\in I}S_i$. \\
      If $M\in \mathbb{F}_\sigma$ and $K\lneq M$, then $\sigma (S_i)=0$ for each $i\in I$ and $M/K\cong \bigoplus\limits _{i\in J}S_i$ for some $J\subseteq I$, so $\sigma (M/K)=0\neq M/K$. Therefore $M\in \mathscr{P}_\sigma $.\\
      On the othe hand, if $M\in \mathscr{P}_\sigma$, then  for each $i\in I$ there exists $K\lneq M$ such that $M/K\cong S_i$. So $\sigma(S_i)\neq S_i$. Then $\sigma(S_i)=0$ for each $i\in I$, so $\sigma(M)=0$.\\
      Therefore $\mathscr{P}_\sigma = \mathbb{F}_\sigma$.
    \end{itemize}
\end{proof}

\noindent\textbf{Luis Fernando Garc\'ia-Mora}\\
Departamento de Matemáticas, Facultad de Ciencias,  \\ 
Universidad Nacional Autónoma de México, Mexico City, Mexico. \\ 
\textbf{e-mail:} \textit{lu1sgarc1agm1995@gmail.com}

\noindent\textbf{Hugo Alberto Rinc\'on-Mej\'ia}\\
Departamento de Matemáticas, Facultad de Ciencias,  \\ 
Universidad Nacional Autónoma de México, Mexico City, Mexico. \\ 
\textbf{e-mail:} \textit{hurincon@gmail.com}\\


\begin{thebibliography}{X}



\bibitem{OSLOMC} \textsc{A. Alvarado-Garc\'ia and  H. Rinc\'on-Mej\'ia},  \textit{On Some Lattices of Module Classes}, 
J. Algebra Appl. 05:01, 105-117 2006.





\bibitem{OTLONACC} \textsc{A. Alvarado-Garc\'ia, H. Rinc\'on-Mej\'ia and J. R\'ios-Montes},  \textit{On the Lattices of Natural and Conatural Classes in R-Mod}, 
Comm. Algebra, 29:2, 541–556, 2001.


\bibitem{An} \textsc{F. W. Anderson} and \textsc{K. R. Fuller}, 
\textit{ Rings and Categories of Modules}, Second Edition, 
Springer-Verlag,  New York, 1992.





\bibitem{Bic} \textsc{L. Bican , T. Kepka,  and P. Nemec},  \textit{Rings, Modules, and Preradicals}, Lecture Notes in Pure and Applied Mathematics. New York: M. Dekker, 1982.

\bibitem{Bic1} \textsc{L. Bican , P. Jambor, T. Kepka,  and P. Nemec}, \textit{Prime and coprime modules}, Fund. Math., 107.1: 33-45, 1980.

\bibitem{Smi} \textsc{S. Çeken, M. Alkan and P. F. Smith}, \textit{Second Modules Over Noncommutative Rings}, 
 Comm. Algebra, 41:1, 83-98, 2013.


\bibitem{Clark} \textsc{J. Clark, C. Lomp, N. Vanaja and R. Wisbauer,} \textit{Lifting Modules: Supplements and Projectivity in Module Theory}, Front. Math., Birkhauser, 2006. 

\bibitem{Da}\textsc{J. Dauns,  and Y. Zhou}, \textit{Classes of Modules}, Chapman Hall/CRC, 2006.

\bibitem{pre} \textsc{R. Fern\'andez-Alonso, F. Raggi, J. R\'ios, H. Rinc\'on and C. Signoret}, \textit{The lattice structure of preradicals} 
Comm. Algebra, 30(3), 1533–1544, 2002.

\bibitem{preII} \textsc{R. Fern\'andez-Alonso, F. Raggi, J. R\'ios, H. Rinc\'on and C. Signoret}, \textit{The lattice structure of preradicals II. Partitions}, 
J. Algebra Appl. 1(2), 201–214, 2002.

\bibitem{preIII} \textsc{R. Fern\'andez-Alonso, F. Raggi, J. R\'ios, H. Rinc\'on and C. Signoret}, \textit{The lattice structure of preradicals III. Operators}, 
 J. Pure Appl. Algebra,  190, 251–265, 2004.



\bibitem{Sec} \textsc{L. F. Garc\'ia-Mora and H. A. Rinc\'on Mej\'ia}, \textit{Second modules relative to subclasses of preradicals of $R$-Mod}, Int. Electron. J. Algebra, Vol. 36, 2024.DOI: 10.24330/ieja.1476650












\bibitem{primepre} \textsc{F. Raggi, J. R\'ios, H. Rinc\'on, R. Fern\'andez-Alonso and C. Signoret}, \textit{Prime and irreducible preradicals}, 
J. Algebra Appl., Vol 4, 2011. 





\bibitem{Es} \textsc{B. Stenström}, \textit{Rings of Quotients: An Introduction to Methods of Ring Theory}. Grundlehren Der Mathematischen Wissenschaften. Berlin, Alemania, Springer-Verlag, 1975.



  \bibitem{Yas} S. Yassemi, \textit{The dual notion of prime submodules}, Arch. Math. (Brno), 37(2001), 273–278.

\end{thebibliography}
\end{document}